\documentclass{article}
\usepackage{texdraw}

\newtheorem{thm}{Theorem}[section]

\newtheorem{rem}[thm]{Remark}

\newcommand{\be}{\begin{equation}}
\newcommand{\ee}{\end{equation}}
\newcommand{\ben}{\begin{enumerate}}
\newcommand{\een}{\end{enumerate}}
\newcommand{\pa}{{\partial}}

\newcommand{\R}{{\rm R}}
\newcommand{\e}{{\epsilon}}

\newcommand{\pxi}{ {\pa \over \pa x^i}}

\font\BBb=msbm10 at 12pt
\newcommand{\Bbb}[1]{\mbox{\BBb #1}}

\title{\Large Two-Dimensional Finsler Metrics\\
 with  Constant Curvature}
\author{Zhongmin Shen}

\date{July,  2001}

\begin{document}

\maketitle

\begin{abstract}
We construct infinitely many two-dimensional Finsler metrics on $\Bbb S^2$ and $\Bbb D^2$ with non-zero constant flag curvature. They are all not locally projectively flat.
\end{abstract}

\section{Introduction}

In Finsler geometry, the flag curvature is an analogue of sectional curvature in Riemannian geometry.  A natural problem is to study and characterize Finsler metrics of constant (flag) curvature. There are  only three local Riemannian metrics of constant curvature, up to a scaling. However there are lots of non-Riemannian Finsler metrics of constant curvature, due to the 
non-Riemannian features of general Finsler metrics. The first set of non-Riemannian Finsler metrics of constant curvature are 
the Hilbert-Klein metric and the Funk metric on a strongly convex domain. The Funk metric is  positively complete and non-reversible with ${\bf K}=-1/4$ and the Hilbert-Klein metric is  complete and reversible with ${\bf K}=-1$.
Both metrics are locally projectively flat \cite{Ok}\cite{Sh1}. 
P. Funk first completely determined the local structure of  two-dimensional projectively flat Finsler metrics with constant curvature
 \cite{Fk1}\cite{Fk2}. 
R. Bryant showed that up to diffeomorphism, there is exactly  a 2-parameter family of locally projectively flat Finsler metrics on $\Bbb S^2$ with ${\bf K}=1$  and the only reversible one is the standard Riemannian metric \cite{Br1}\cite{Br2}. Later on,  he extended his construction to higher dimensional spheres $\Bbb S^n$ \cite{Br3}. 
Recently, the author has completely determined the local  structure of projectively flat analytic Finsler metrics of constant curvature in higher dimensions \cite{Sh4}. Our method is different from Funk's.

The next problem is to classify non-projectively flat
Finsler metrics of constant curvature. This problem turns out to be very 
difficult. The very first step might be  to construct as many examples as possible.
In 2000, D. Bao and the author first constructed a family of non-projectively flat Finsler metrics on $\Bbb S^3$ with ${\bf K}=1$ using the Lie group structure of ${\Bbb  S}^3$ \cite{BaSh}. 
Our examples are in the form $F= \alpha+\beta$, where $\alpha ({\bf y})
=\sqrt{a_{ij}(x)y^iy^j }$ is a Riemannian metric and $\beta({\bf y})
= b_i(x) y^i$ is a  1-form. Finsler metrics in this form are called Randers metrics \cite{Ra}.  
Recently, the author has just constructed 
an incomplete non-projectively flat  Randers metric with ${\bf K}=0$ in each dimension\cite{Sh3}. 

The main technique in \cite{Sh3} is described as follows. Given a Finsler metric 
$\Phi$  and a vector field ${\bf v}$ on a manifold $M$, define a function
$F: TM \to [0, \infty)$ by
\be
 \Phi \Big ( { {\bf y}\over F({\bf y})} -\e{\bf v}_p \Big ) = 1, \ \ \ \ \ \ {\bf y}\in T_pM.\label{eq}
\ee
where $\e $ is a constant.
$F$ is a Finsler metric when $\e$ is small.
An important relationship between $\Phi$ and $F$ is that their (Busemann-Hausdorff) volume forms are equal, $dV_{\Phi}=dV_F$ \cite{Sh3}.
By choosing an appropriate Finsler metric $\Phi$ and an appropriate  vector field 
${\bf v}$, one  obtains a Finsler metric $F$ with the same  curvature properties as $\Phi$. 

In this paper, we are going to employ this technique to construct a family 
of Randers metrics on $\Bbb S^2$ with ${\bf K}=1$ and a family of Randers metrics on 
a disk $\Bbb D^2(\rho)$ with ${\bf K}=-1$ or ${\bf K}=-1/4$. They  are all not locally projectively flat. These examples show that the classification problem of non-projectively flat  Finsler metrics 
of constant curvature is  very difficult.

Now let us describe our examples.
Let $\Phi({\bf y}) = \sqrt{ h( {\bf y}, {\bf y}) }$ denote the
standard Riemannian metric on the unit sphere $\Bbb S^2$ and ${\bf v}$ denote the vector field on $\Bbb S^2$ defined by
\be
 {\bf v}_p = \Big (-y,\;  x, \; 0\Big ) \ \ \ \ {\rm at}  \ p = \Big (x, \; y, \; z\Big )\in \Bbb S^2,\label{external1}
\ee
Define $F: T{\Bbb S^2} \to [0, \infty)$ by (\ref{eq}). Then 
 $F= \alpha +\beta $ is a Randers metric, where $\alpha=\alpha({\bf y})$ and $\beta =\beta({\bf y})$ are given by
\be
\alpha: ={ \sqrt{ \e^2 h( {\bf v}, {\bf y} )^2 + 
h( {\bf y}, {\bf y} ) \Big ( 1- \e^2 h( {\bf v}, {\bf v} ) \Big ) }
\over 1 -\e^2 h( {\bf v}, {\bf v} ) }, \ \ \ \ 
\beta: =  - {\e h( {\bf v}, {\bf y} )\over 
1-\e^2h( {\bf v}, {\bf v} ) }. \label{Spherical}
\ee
$F$ is defined on the whole sphere for $|\e | < 1$
and it is defined only on the open  disks  around the north pole and south
pole  with radius $\rho =\sin^{-1}(1/|\e|) $
for $|\e | \geq 1$. Note that when $\e =0$, $F=\Phi$ is the standard Riemannian metric on $\Bbb S^2$.

\begin{thm}\label{thm1}
Let $F=\alpha+\beta$ be any Finsler metric on $\Bbb S^2$ given in (\ref{Spherical}).
It 
 has the following properties
\ben
\item[(a)] ${\bf K}= 1$;
\item[(b)] ${\bf S}=0$;
\item[(c)] $ F$ is not locally projectively flat unless $\e =0$;
\item[(d)] the Gauss curvature $\bar{\bf K}$ of $\alpha$
is not a constant  unless $\e=0, \pm 1$. When $\e=0$, 
$\bar{\bf K}=1$; When $\e=\pm 1$, 
$\bar{\bf K}= -4 .$
\een
\end{thm}

According to Yasuda-Shimada\cite{YaSh}, if a Randers metric 
$F=\alpha+\beta$ is of positive constant curvature , then $\beta$ must be a Killing form of constant length with respect to $\alpha$. However,  $\beta$ in (\ref{Spherical}) is not a Killing form when $\e \not=0$. Thus Theorem \ref{thm1} is inconsistent with  Yasuda-Shimada's result in dimension two.

\bigskip

Similarly, 
let $\Phi ({\bf y}) =\sqrt{ h( {\bf y}, {\bf y} ) } $ denote the standard Klein metric on
the unit disk $\Bbb D^2$ and ${\bf v}$ denote the vector field on
$\Bbb D^2$ defined by
\be
 {\bf v}_p = (-y,\ x) \ \ \ \ {\rm at} \ p=(x, y)\in \Bbb D^2.\label{external2}
\ee
Define $F: T{\Bbb D}^2 \to [0, \infty)$ by (\ref{eq}). Then  
$F= \alpha+\beta$ is a Randers metric, where $\alpha=\alpha({\bf y})$ and $\beta =\beta({\bf y})$ are given by
\be
\alpha: = { \sqrt{ \e^2 h( {\bf v}, {\bf y} )^2 + 
h( {\bf y}, {\bf y} ) \Big ( 1- \e^2 h( {\bf v}, {\bf v} ) \Big ) }
\over 1 -\e^2 h( {\bf v}, {\bf v} ) }, \ \ \ \ \ 
\beta: = - {\e h( {\bf v}, {\bf y} )\over 
1-\e^2h( {\bf v}, {\bf v} ) }.\label{Hyperbolic}
\ee 
$F$ is a Finsler metric defined on the disk $\Bbb D^2(\rho)$ 
with radius $\rho=1/\sqrt{1+\e^2} $. Note that when $\e=0$,
$F$ is the Klein metric on the unit disk.

\begin{thm}\label{thm2} Let $F=\alpha+\beta$ be the  Finsler metric on 
the disk $\Bbb D^2(\rho)$ given in (\ref{Hyperbolic}). 
It has the following properties
\ben
\item[(a)] ${\bf K}= - 1$;
\item[(b)] ${\bf S}=0$;
\item[(c)] $ F$ is not locally projectively if $\e \not=0$;
\item[(d)] the Gauss curvature $\bar{\bf K}$ of $\alpha$ is not constant
unless $\e=0$. When $\e=0$, $\bar{\bf K}=-1$.
\een
\end{thm}

\bigskip
According to Yasuda-Shimada \cite{YaSh} if a Randers metric $F=\alpha+\beta$ 
is of negative constant curvature, then the Riemannian metric $\alpha$ is of negative constant curvature.  However, the Randers metric defined in (\ref{Hyperbolic}) do not have this property when $\e\not=0$. Thus
Theorem \ref{thm2} is inconsistent  with  Yasuda-Shimada's result  in dimension two.

\bigskip

Besides the Klein metric, the hyperbolic metric can be expressed in many other forms, such as the Poincare metric and the one arising from the proof of Theorem \ref{thm1}. We can use them  to construct many  non-projectively flat
 Finsler metrics with negative  constant curvature. See Remark \ref{remK=-1} below.

\bigskip
Finally, let $\Phi ({\bf y}) 
= \sqrt{h( {\bf y}, {\bf y}) } + h( {\bf u}, {\bf y}) $ denote the Funk metric on the unit disk $\Bbb D^2$,
where $h $ is the Klein metric on $\Bbb D^2$ and ${\bf u}= (1-x^2-y^2)(x {\pa \over \pa x} + y {\pa \over \pa y} )\in T_{(x, y)} \Bbb D^2$ is a vector field.
$\Phi({\bf y})$, ${\bf y}\in T_p\Bbb D^2$,  is defined by 
\be
 {{\bf y}\over \Phi({\bf y})} + p \in \pa \Bbb D^2.\label{Funkmetric}
\ee 
Let ${\bf v}$ denote the vector field on $\Bbb D^2$ defined by
(\ref{external2}). 
Define $F: T\Bbb D^2\to [0,\infty)$ by (\ref{eq}), i.e.,
\be
\sqrt{  h\Big (  {{\bf y}\over F ({\bf y})} -\e {\bf v}, \; {{\bf y}\over F ({\bf y})} -\e {\bf v}\Big  ) }
+ h\Big ( {\bf u},\;  {{\bf y}\over F ({\bf y})} -\e {\bf v}\Big  ) =1.\label{FunkFunk}
\ee
Then $F=\alpha+\beta$ is a Randers metric on $\Bbb D^2(\rho)$ where $\rho = 1/\sqrt{1+\e^2}$.

\begin{thm}\label{thm3} Let $F=\alpha+\beta$ be the Randers metric on $\Bbb D^2(\rho)$
defined in (\ref{FunkFunk}). It has the following properties:
\ben
\item[(a)] ${\bf K}= -1/4$;
\item[(b)] ${\bf S}= {3\over 2} F$;
\item[(c)] $F$ is not locally projectively flat unless $\e =0$;
\item[(d)] the Gauss curvature $\bar{\bf K}$ of $\alpha$ is not a constant unless $\e =0$. When $\e =0$, $\bar{\bf K}= -1$.
\een
\end{thm}

Again Theorem \ref{thm3} is inconsistent with Yasuda-Shimada's result in dimension two, since $\alpha$ does not have constant curvature when $\e \not=0$.

\bigskip

In a recent paper by Bao-Robles \cite{BaRo}, they characterize Randers metrics with constant curvature by three equations. Moreover, they use the technique in \cite{Sh3} to construct two-dimensional Randers metrics with the Gauss curvature ${\bf K}= {\bf K}(x)$ independent of the directions, a three-dimensional Randers metric on ${\rm S}^3$ with ${\bf K}=1$
and a three-dimensional Randers metric on ${\rm B}^3$ with ${\bf K}=-1$. 
We should point out that their example on ${\rm S}^3$ is not equivalent to that in \cite{BaSh}.
With the examples in \cite{BaSh}\cite{BaRo}\cite{Sh3}, now we have 
non-projectively flat Randers metrics with constant curvature of any sign in higher dimensions.

\section{Preliminaries}\label{sectionpre}

Let $F$ be a Finsler metric on a manifold $M$. 
In a standard local coordinate system $(x^i, y^i)$ in $TM$, 
$F= F(x, y)$ is a function of $(x^i, y^i)$. Let
\[ g_{ij} = {1\over 2} [F^2]_{y^iy^j}\]
and $(g^{ij}):=(g_{ij})^{-1}$.
The geodesics of $F$ are characterized locally by
\[ {d^2 x^i\over dt^2} + 2 G^i \Big ( x, {dx\over dt} \Big )=0,
\]
where
\[ G^i = {1\over 4} g^{ik} \Big \{ 2 {\pa g_{pk}\over \pa x^q} 
- {\pa g_{pq}\over \pa x^k} \Big \} y^py^q.\]
The coefficients of the Riemann curvature ${\bf R}_{\bf y}= R^i_{\ k} dx^k 
\otimes \pxi$ are given by
\be
R^i_{\ k} = 2 {\pa G^i \over \pa x^k}
- y^j {\pa^2 G^i\over \pa x^j \pa y^k} + 2 G^j {\pa G^i\over \pa y^j \pa y^k} - {\pa G^i\over \pa y^j} {\pa G^j \over \pa y^k}. \label{Rik}
\ee
$F$ is said to be of constant curvature  ${\bf K}=\lambda$, if 
\[ R^i_{\ k} = \lambda \Big \{ F^2 \delta^i_k - F F_{y^k} y^i \Big \}.\]

\bigskip
When $F=\sqrt{a_{ij}(x)y^iy^j}$ is a Riemannian metric,
$R^i_{\ k}= R^{\ i}_{j \ kl}(x) y^j y^l$, where $R^{\ i}_{j\ kl}(x)$ denote the coefficients of the usual Riemannian curvature tensor.
Thus the quantity ${\bf R}_{\bf y}$ in Finsler geometry is  still called the Riemann curvature.

There are many interesting non-Riemannian quantities in Finsler geometry.
In this paper, we will only discuss the  S-curvature \cite{Sh1}.
Express the (Busemann-Hausdorff) volume form of $F$ by
\[ dV_F = \sigma(x) dx^1 \cdots dx^n.\]
The S-curvature is defined by
\be
{\bf S}({\bf y}) := {\pa G^i \over \pa y^i} (x,y)- { y^i \over \sigma(x) } {\pa \sigma \over \pa x^i}(x).\label{localS}
\ee
See \cite{Sh1} for a related discussion on the S-curvature. 

\bigskip
Randers metrics are among the simplest non-Riemannian Finsler metrics, so that many well-known geometric quantities are computable.

Let  $F = \alpha+\beta$ be a Randers metric on a manifold $M$,
where 
\[\alpha (y)=\sqrt{a_{ij}(x) y^iy^j}, \ \ \ \ \ \beta(y)= b_i(x) y^i\]
with 
$\|\beta\|_x :=\sup_{y \in T_xM} \beta(y)/\alpha(y) < 1$. 
Define $b_{i|j}$ by
\[ b_{i|j} \theta^j := db_i -b_j \theta_i^{\ j},\]
where 
$\theta^i :=dx^i$ and 
$\theta_i^{\ j} :=\tilde{\Gamma}^j_{ik} dx^k$ denote the Levi-Civita  connection forms of $\alpha$.
Let
\[ r_{ij} :={1\over 2} \Big ( b_{i|j}+b_{j|i}\Big ) , \ \ \ \ \ 
s_{ij}:= {1\over 2} \Big (b_{i|j}-b_{j|i} \Big ),\]
\[s^i_{\ j}:= a^{ih}s_{hj}, \ \ \ \ \  s_j:=b_i s^i_{\ j}, \ \ \ \ \   e_{ij} := r_{ij}
+ b_i s_j + b_j s_i.\]
Then $G^i$  are given by
\be
 G^i = \bar{G}^i + {e_{00} \over 2F} y^i -s_0y^i+ \alpha s^i_{\ 0},\label{Gi}
\ee
where $e_{00}:= e_{ij}y^iy^j$, $s_0:=s_iy^i$, $s^i_{\ 0}:=s^i_{\ j} y^j$
and $\bar{G}^i$ denote the geodesic coefficients of $\alpha$.
See \cite{AIM}. 

According to Lemma 3.1 in \cite{ChSh}, 
\be
{\bf S}= c (n+1) F  \ \ \ \ \Longleftrightarrow \ \ \ \ e_{00} = 2c (\alpha^2-\beta^2). \label{Se}
\ee
where $c = c(x)$ is a scalar function. See also Proposition 5.1 in \cite{Sh3} in the case
when $c=0$.

Assume that  
${\bf S} = c (n+1) F$ for some constant $c$. 
Then
\be
G^i = \bar{G}^i + c (\alpha-\beta) y^i - s_0 y^i + \alpha s^i_{\ 0}.\label{Gi2}
\ee
By a direct computation, one obtains a formula for  the Riemann curvature 
is given by
\begin{eqnarray}
R^i_{\ k} &=& \bar{R}^i_{\ k}+3c^2\Big (\alpha^2 \delta^i_k- y^i y_k\Big )  -c^2\beta \Big ( \beta \delta^i_k -b_ky^i \Big )\nonumber\\
&& 
 +\Big ( s_{0|0} \delta^i_k -s_{0|k}  y^i\Big )  + s_0\Big ( s_0 \delta^i_k -s_k y^i \Big )  +\Big ( s_{k|0} - s_{0|k}  \Big ) y^i \nonumber \\
&& -\Big (  \alpha^2s^i_{\ j}s^j_{\ k} -y_k s^i_{\ j}s^j_{\ 0} \Big )   +6 c s_{k0} y^i +3 s_{k0} s^i_{\ 0}\nonumber\\
&& - \Big \{(c^2 \beta + 2 cs_0 +s_js^j_{\ 0}) \Big ( \alpha^2 \delta^i_k -y_k y^i\Big )  + c^2 \alpha^2 \Big ( \beta \delta^i_k -b_k y^i \Big )\nonumber\\
&&  +2c \alpha^2 \Big ( s_0 \delta^i_k -s_k y^i \Big )  - \Big ( \alpha^2 s^i_{\ 0|k} - y_k s^i_{\ 0|0} \Big )\nonumber\\
&& + \alpha^2 \Big ( s_j s^j_{\ 0} \delta^i_k - s_j s^j_{\ k} y^i \Big ) +\alpha^2 \Big ( s^i_{\ k|0} - s^i_{\ 0|k} \Big )   \Big \}\alpha^{-1}.
\label{localRik}
\end{eqnarray}

\bigskip
Taking the trace of $R^i_{\ k}$, we  obtain a formula for the Ricci curvature ${\bf Ric}$ of $F$ which is expressed in terms of the Ricci curvature $\overline{\bf Ric}$ of $\alpha$ and the covariant derivatives of $\beta$ with respect to $\alpha$. 
\begin{eqnarray}
{\bf Ric} & = & \overline{\bf Ric} 
+ (n-1) \Big \{ c^2 (\alpha^2+\beta^2) + 2c^2 (\alpha^2-\beta^2)  + s_{0|0} + s_0s_0 \Big \}\nonumber \\
&& + 2 s_{k0}s^k_{\ 0} -\alpha^2 s^k_{\ j} s^j_{\ k} \nonumber\\
&&  + \Big \{ 2 s^k_{\ 0|k} - (n-1) \Big ( 4c s_0 + 2 s_js^j_{\ 0}+   2c^2 \beta\Big )   \Big \}\alpha. \label{RikRik}
\end{eqnarray}

\section{Proof of Theorem \ref{thm1}}

The Finsler metric in Theorem \ref{thm1} is constructed 
by solving  the equation (\ref{eq}), i.e.,
\be 
 \Phi \Big ({{\bf y}\over F({\bf y})} - \e {\bf v}\Big ) = \sqrt{h\Big ( { {\bf y}\over F({\bf y})} - \e {\bf v} ,\;{ {\bf y}\over F({\bf y})} - \e {\bf v} \Big ) }= 1.\label{eq1}
\ee

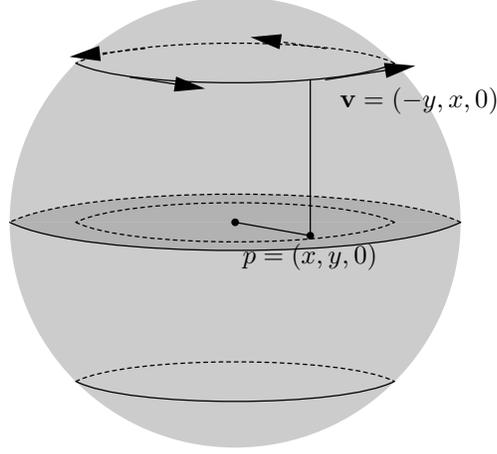
\begin{center}
\begin{texdraw}
\drawdim cm
\linewd  0.02

\move(0 0)
\fcir f:0.8 r:3
\move(-3 0)
\clvec(-2 -0.5)(2  -0.5)(3 0)
\lfill f:0.7
\lpatt(0.05 0.05)
\clvec(2 0.5)(-2 0.5)(-3 0)
\lfill f:0.7
\move(0 0)
\fcir f:0.0 r:0.05

\lpatt()
\move(-2.12 2.12)
\clvec(-1.42 1.77)(1.42 1.77)(2.12 2.12)
\lpatt(0.05 0.05)
\clvec(1.42 2.47)(-1.42 2.47)(-2.12 2.12)

\lpatt(0.05 0.05)
\move(-2.12 -2.12)
\clvec(-1.42 -1.77)(1.42 -1.77)(2.12 -2.12)
\lpatt()
\clvec(1.42 -2.47)(-1.42 -2.47)(-2.12 -2.12)

\lpatt(0.05 0.05)
\move(-2.12 0)
\clvec(-1.42 -0.35)(1.42 -0.35)(2.12 0)
\lpatt(0.05 0.05)
\clvec(1.42 0.35)(-1.42 0.35)(-2.12 0)

\lpatt(0.05 0.05)
\move(1.2 2.315)
\arrowheadtype t:F \avec(0.2 2.46)

\move(-1.2 2.315)
\arrowheadtype t:F \avec(-2.2 2.2)

\lpatt()
\move(1.2 1.9)
\arrowheadtype t:F \avec(2.4 2.08)
\move(-1.4 1.93)
\arrowheadtype t:F \avec(-0.4 1.78)

\move(1.4 1.8)
\textref h:L v:T \htext{${\bf v}=(-y, x, 0)$}

\move(0 0)
\lvec(1 -0.18)
\fcir f:0.0 r:0.05
\lvec(1 1.9)
\move(1 -0.28)
\textref h:C v:T \htext{$p=(x, y, 0)$}

\move(0 3.5)
\fcir f:1 r:0.0
\move(0 -0.5)
\fcir f:1 r:0.0
\end{texdraw}
\end{center}

Let $\psi: \R^2 \to \Bbb S^2_{+}$ by
\[ \psi (x, y) := \Big ( {x\over \sqrt{ 1+x^2+y^2}},\ {y\over \sqrt{ 1+ x^2+y^2}},\ { 1\over \sqrt{1+x^2+y^2} }\Big ).\]
With this map, the standard Riemannian metric $\Phi$ on 
$\Bbb S^2$ can be expressed on $\R^2$ by
\[ \Phi({\bf y}) = {\sqrt{ (u^2+v^2) + (xv-yu)^2  }\over 1+x^2 +y^2},\]
where ${\bf y}= u{\pa \over \pa x} + v {\pa \over \pa y} \in T_{(x, y)}\R^2$. The Finsler metric  defined by (\ref{eq1}) is a 
 Randers metric $F=\alpha+\beta$, where 
$\alpha=\alpha({\bf y})$ and $\beta =\beta({\bf y})$ are given by
 \begin{eqnarray*}
{\alpha} : & = & {\sqrt{ \Big ( 1+(1-\e^2)(x^2+y^2)\Big )(u^2+v^2)
 + \Big ( 1+\e^2 + x^2+y^2\Big ) (xv-yu)^2  }
\over \Big ( 1+(1-\e^2)(x^2+y^2) \Big ) \sqrt{ 1+x^2+y^2} }\\
{\beta} : & = & - { \e (xv-yu)\over 1+(1-\e^2) (x^2+y^2) }.
\end{eqnarray*}
Note that when $|\e|>1$,  $F$ is defined only on the open disk $\Bbb D^2(r)$ of radius 
$r = 1/\sqrt{\e^2-1}$. The corresponding domain  on $\Bbb S^2$ 
is a metric disk $B(\rho)$ around the north pole with radius $\rho= \sin^{-1} (1/|\e|)$.

To compute the curvatures of $F$, we express it in a polar coordinate system,
$ x = r \cos(\theta), \ y = r \sin(\theta).$
For ${\bf y}= \mu {\pa \over \pa r} + \nu {\pa \over \pa \theta}$,
$\alpha=\alpha({\bf y})$ and $\beta =\beta({\bf y})$ are given by
\begin{eqnarray*}
\alpha & = & {\sqrt{\Big (1+(1-\e^2) r^2\Big )\mu^2+ r^2\Big (1+r^2\Big )^2\nu^2}\over 
\Big (1+(1-\e^2)r^2\Big )\sqrt{1+r^2}}\\
\beta& = & - { \e r^2 \nu\over 1+ (1-\e^2)r^2}
\end{eqnarray*}
Express $\alpha =\sqrt{ a_{11}\mu^2+ a_{12} \mu\nu + a_{21} \nu \mu+ a_{22}\nu^2}$
and $\beta= b_1 \mu +b_2 \nu$, where
\[ a_{11}= {1\over (1+r^2)(1+(1-\e^2)r^2)}, 
\ \ \ \ a_{12}=0 = a_{21}, \ \ \ \ a_{22} = { r^2 (1+r^2)\over (1+(1-\e^2)r^2)^2},\]
\[ b_1 =0, \ \ \ \ \  b_2 = - {\e r^2 \over 1 + (1-\e^2) r^2 }.\]
The geodesic coefficients $\bar{G}^1 $ and $\bar{G}^2$  of $\alpha$ are given by
\begin{eqnarray*}
\bar{G}^1 & = &{ \Big (1 +(1  -\e^2)(1 + 2  r^2  ) \Big )r \over 2(1+r^2)\Big (1+(1-\e^2)r^2\Big )}\; \mu^2  - { (1+r^2)\Big (1+2 r^2-(1-\e^2)r^2\Big  ) r \over 2\Big (1+(1-\e^2)r^2\Big )^2}\;\nu^2\\
\bar{G}^2 & = &{ 1+2 r^2- (1-\e^2) r^2 \over (1+r^2)\Big (1+(1-\e^2)r^2\Big ) r }\; \mu\nu
\end{eqnarray*}
We immediately obtain the Gauss curvature $\bar{\bf K}$ of $\alpha$,
\[ \bar{\bf K} ={1-5 \e^2 + (1-\e^4) r^2 \over 1+ (1-\e^2)r^2 } .\]
Note that for $\e = \pm 1$, $\alpha$ has negative constant Gauss curvature,
\[ \bar{\bf K}= - 4.\]

Now we are going to find the geodesic coefficients $G^1$ and $G^2$ of $F$.
By (\ref{Gi}), we first compute $r_{ij}, s^i_{\ j}$ and $s_i $, etc.
A direct computation yields that
\begin{eqnarray*}
r_{11} & = & 0 \ \ = \ \ r_{22}\\
r_{12} & = & { \e^3 r^3\over (1+r^2)\Big (1+ (1-\e^2)r^2\Big )^2} \ \ = \ \ r_{21}\\
s_{11} & = & 0 \ \ = \ \ s_{22}\\
s_{12} & = & { \e r \over \Big ( 1+ (1-\e^2)r^2\Big  )^2 } = - s_{21}\\
s^1_{\ 1} & = & 0 \ \ = \ \ s^2_{\ 2}\\
s^1_{\ 2} & = & { \e r (1+r^2)\over 1 + (1-\e^2)r^2 }\\
s^2_{\ 1} & = & - {\e \over r (1+r^2) }\\ 
s_1 & = & { \e^2 r\over (1+r^2) \Big ( 1+ (1-\e^2)r^2\Big ) }\\
s_2 & = & 0.
\end{eqnarray*}
We obtain that
\[ e_{ij} := r_{ij} + b_i s_j + b_j s_i= 0\]
This is equivalent to that ${\bf S}=0$.  By (\ref{Gi2}) and the above identities, we obtain
\begin{eqnarray*}
G^1 & = & \bar{G}^1   -{\e^2 r  \over (1+r^2)\Big (1+(1-\e^2)r^2\Big )}\; \mu^2 + {\e r (1+r^2)\over 1+(1-\e^2)r^2 } \; \alpha \; \nu \\
G^2 & = &\bar{G}^2 -{\e^2 r  \over (1+r^2)(1+(1-\e^2)r^2)}\; \mu \nu
 - {\e \over r (1+r^2) }\; \alpha\; \mu
\end{eqnarray*}
Plugging them into (\ref{Rik}), we obtain
\be
 R^i_{\ k} = F^2 \Big \{ \delta^i_k -  {F_{y^k}\over F}  y^i\Big \}.\label{Rik=1}
\ee
We conclude that  the Gauss curvature ${\bf K}=1$.

We can also use (\ref{RikRik}) and the above identities to verify 
that ${\bf K}=1$.
To do so, it suffices to compute $s_{0|0}$ and $s^k_{\ 0|k}$. They are given by
\begin{eqnarray*}
s_{0|0} & = & { \e^2 \Big ( 1 - (1-\e^2)r^4 \Big ) \over 
\Big ( 1+r^2 \Big )^2 \Big ( 1 + (1-\e^2) r^2 \Big )^2 }\mu^2
+ { \e^2 r^2 \Big ( 1+ (1+\e^2)r^2 \Big ) \over \Big ( 1 + (1-\e^2)r^2 \Big )^3 }\nu^2 \\
s^k_{\ 0|k} & = & - {\e (1-\e^2) r^2 ( 1+r^2) \nu \over \Big ( 1+(1-\e^2)r^2 \Big )^2 }.
\end{eqnarray*}
Plugging them into (\ref{RikRik}) gives 
\[ {\bf Ric} =  F^2 .\]
We conclude that ${\bf K}={\bf Ric}/F^2 = 1$.

\bigskip

\begin{rem}\label{remK=1}
{\rm
Express the spherical metric in a radial form 
\[ \Phi ({\bf y}) 
= \sqrt{ u^2 + \sin^2(r) v^2 },
\]
where ${\bf y}= u {\pa \over \pa r} + v {\pa \over \pa \theta}\in T_{(r, \theta)}((0, \infty)\times {\rm S}^1 )$.
 Take 
${\bf v}= {\pa \over \pa \theta}\in T_{(r, \theta)}((0, \infty)\times {\rm S}^1 )$ and define $F$ by (\ref{eq}). We obtain
\be
F = {\sqrt{ \Big ( 1-\e^2 \sin^2(r) \Big ) u^2 + \sin^2(r) v^2 } -\e \sin^2(r) v\over 1-\e^2 \sin^2(r) }.
\ee
$F$ satisfies that ${\bf K}=1$ and ${\bf S}=0$, but it is not locally projectively flat.
}
\end{rem}

\section{Proof of Theorem \ref{thm2}}

The Finsler metric in Theorem \ref{thm2} is also constructed 
by solving the equation (\ref{eq}), i.e.,
\be 
 \Phi \Big ({{\bf y}\over F({\bf y})} - \e {\bf v}\Big ) = \sqrt{h\Big ( { {\bf y}\over F({\bf y})} - \e {\bf v} ,\;{ {\bf y}\over F({\bf y})} - \e {\bf v} \Big ) }= 1.\label{eq2}
\ee

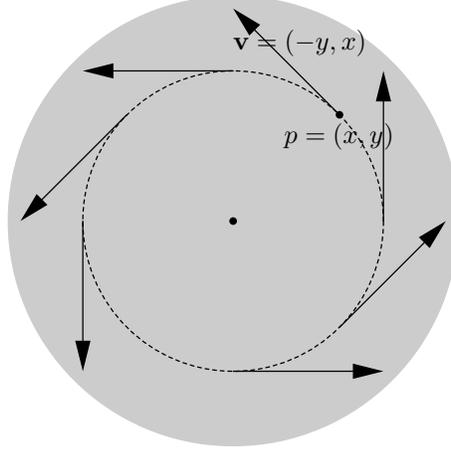
\begin{center}
\begin{texdraw}
\drawdim cm
\linewd  0.02

\move(0 0)
\fcir f:0.8 r:3
\fcir f:0.0 r:0.05

\lpatt(0.05 0.05)
\lcir r:2

\lpatt()
\move(2 0)
\arrowheadtype t:F \avec(2 2)
\move(0 2)
\arrowheadtype t:F \avec(-2 2)
\move(-2 0)
\arrowheadtype t:F \avec(-2 -2)
\move(0 -2)
\arrowheadtype t:F \avec(2 -2)
\move(1.414 1.414)
\arrowheadtype t:F \avec(0 2.828)
\move(-1.414 1.414)
\arrowheadtype t:F \avec(-2.828 0)
\move(1.414 -1.414)
\arrowheadtype t:F \avec(2.828 0)
\move(1.414 1.414)
\fcir f:0.0 r:0.05
\move(1.414 1.3)
\textref h:C v:T \htext{$p=(x, y)$}
\move(0 2.2)
\textref h:L v:B \htext{${\bf v}= (-y, x)$}

\end{texdraw}
\end{center}

The Klein metric $\Phi$  on $\Bbb D^2$ is given by
\[ \Phi ({\bf y}) = {\sqrt{  (u^2+v^2) -(xv-yu)^2 }\over 1- (x^2+y^2) },\]
where 
${\bf y}=(u, v)\in T_{(x, y)}\R^2$. 
The Finsler  metric defined by 
(\ref{eq2}) is a Randers metric $F=\alpha+\beta$,
where 
$\alpha =\alpha({\bf y})$ and $\beta=\beta({\bf y})$ are given by
\begin{eqnarray*}
{\alpha} : & = & {\sqrt{ \Big(1-(1+\e^2)(x^2+y^2)\Big )(u^2+v^2)
-\Big ( 1-\e^2 -(x^2+y^2) \Big ) (xv-yu)^2 }
\over \Big ( 1- (1+\e^2)(x^2+y^2) \Big ) \sqrt{ 1-x^2-y^2} }\\
{\beta} : & = & - { \e (xv-yu)\over 1- (1+\e^2) (x^2+y^2) }.
\end{eqnarray*}

To compute the curvatures of $F$, we take a polar coordinate system,
$x = r \cos (\theta), y = r\sin(\theta)$.
For a vector ${\bf y} = \mu {\pa \over \pa r} +\nu {\pa \over \pa \theta}$,
$\alpha =\alpha({\bf y})$ and $\beta =\beta({\bf y})$ are given by
\begin{eqnarray*}
\alpha & = & {\sqrt{\Big (1-(1+\e^2)r^2\Big ) \mu^2 + r^2\Big (1-r^2\Big )^2 \nu^2}\over \Big (1-(1+\e^2)r^2\Big )\sqrt{1-r^2}}\\
\beta & = & - {\e r^2 \nu\over 1 - (1+\e^2)r^2 }.
\end{eqnarray*}
Express $\alpha =\sqrt{ a_{11}\mu^2 + a_{12} \mu \nu + a_{21}\nu \mu + a_{22} \nu^2 }$
and $\beta = b_1 \mu + b_2 \nu$, where
\[a_{11} = { 1 \over \Big ( 1-(1+\e^2)r^2 \Big ) (1-r^2)}\
\ \ \ a_{12} =0 = a_{21}, \ \ \ \ a_{22} = { r^2 (1-r^2) \over \Big ( 1-(1+\e^2)r^2 \Big )^2} ,\]
\[ b_1 =0, \ \ \ \ \ b_2 = - {\e r^2 \over 1- (1+\e^2) r^2 }.\]
The geodesic coefficients $\bar{G}^1$ and $\bar{G}^2$ of $\alpha$ are given by
\begin{eqnarray*}
\bar{G}^1 & = & { \Big (1 + (1 +\e^2)(1-2r^2)  \Big ) r \over 
2(1-r^2) \Big (1-(1+\e^2)r^2\Big )}\; \mu^2
- {\Big (1-2r^2+(1+ \e^2) r^2\Big )(1-r^2) r \over 2\Big (1-(1+\e^2)r^2\Big )^2 }\; \nu^2\\
\bar{G}^2 & = & {  1-2r^2 +(1+\e^2) r^2 \over (1-r^2)\Big (1-(1+\e^2)r^2\Big ) r }\; \mu \nu
\end{eqnarray*}
The Gauss curvature $\bar{\bf K}$ of $\alpha$ is given by
\be
 \bar{\bf K} = { -1 -5 \e^2 + (1-\e^4) r^2 \over 1 - (1+\e^2) r^2 }.\label{K1}
\ee
We see that $\bar{\bf K}$ is not a constant unless $\e =0$.

Now we are going to find the geodesic coefficients $G^1$ and $G^2$ of $F=\alpha+\beta$. Let $r_{ij}, s_{ij}, s^i_{\ j}, s_j $ and $e_{ij}$ as above. 
A direct computation yields that 
\begin{eqnarray*}
r_{11} & = & 0 \ \ = \ \ r_{22}\\
r_{12} & = &  {\e^3 r^3\over (1-r^2)\Big (1-(1+\e^2)r^2\Big )^2} \ \ = \ \ r_{21}\\
s_{11} & = & 0 \ \ =\ \  s_{22}\\
s_{12} & = & {\e r \over \Big ( 1 - (1+\e^2)r^2 \Big )^2 } \ \ =  \ \ -s_{21}\\
s^1_{\ 1} & = &  0 \ \ = \ \ s^2_{ \ 2} \\
s^1_{\ 2} & = & { \e r (1-r^2) \over 1- (1+\e^2)r^2 }\\
s^2_{\ 1} & = & - {\e\over (1-r^2) r }\\
s_1 & = & {\e^2 r \over (1-r^2)\Big (1- (1+\e^2) r^2\Big  ) }\\
s_2 & = & 0.
\end{eqnarray*}
We immediately see that 
\[ e_{ij} :=r_{ij}+b_i s_j +b_j s_i=0.\]
Thus the S-curvature vanishes, ${\bf S}=0$.  By (\ref{Gi2}) and  the above identities, we obtain
\begin{eqnarray*}
G^1 & = & \bar{G}^1 - {\e^2 r \over (1-r^2) (1- (1+\e^2)r^2 ) } \; \mu^2 
+  { \e r (1-r^2) \over 1-(1+\e^2)r^2 }\; \alpha\nu\\
G^2 & = & \bar{G}^2 - {\e^2 r  \over (1-r^2) (1- (1+\e^2)r^2 ) }\;\mu\nu
- {\e  \over (1-r^2) r }\; \alpha \mu.
\end{eqnarray*}
Plugging them into (\ref{Rik}), we immediately obtain
\be
R^i_{\ k} = - \Big \{ F^2 \delta^i_k - F F_{y^k} y^i\Big \}.\label{Rik=-1}
\ee
Thus the Gauss curvature ${\bf K}=-1$. 

We can also use (\ref{RikRik}) and the above identities to verify that ${\bf K}=-1$. To do so, it suffices to compute 
$s_{0|0}$ and $s^k_{\ 0|k}$. They are given by
\begin{eqnarray*}
s_{0|0} & = & { \e^2 \Big ( 1 - (1+\e^2) r^4 \Big ) \over \Big ( 1-r^2 \Big )^2 \Big ( 1 - (1+\e^2)r^2 \Big )^2 }\; \mu^2 
+ { \e^2 r^2 \Big ( 1 - (1-\e^2 ) r^2 \Big ) \over \Big (1 - (1+\e^2) r^2 \Big )^3 }\;\nu^2,\\
s^k_{\ 0|k} & = & { \e (1+\e^2) r^2 ( 1 -r^2) \over \Big ( 1 - (1+\e^2)r^2 \Big )^2 } \; \nu.
\end{eqnarray*}
Plugging them into (\ref{RikRik}), we obtain
\[ {\bf Ric} =- F^2.\]
Again, we conclude that ${\bf K}={\bf Ric}/F^2 = -1$.

\bigskip
\begin{rem}\label{remK=-1}
{\rm 
Express the Klein metric in the radial form,
\[ \Phi({\bf y})
= \sqrt{u^2 + \sinh^2 (r) v^2} , 
\]
where ${\bf y}= u {\pa \over \pa r} + v {\pa \over \pa \theta} \in T_{(r, \theta)} ((0, \infty) \times {\rm S}^1 )$.
Take ${\bf v}= {\pa \over \pa \theta} \in T_{(r, \theta)}(\R \times {\rm S}^1 )$ and define $F$ by (\ref{eq}).  We obtain
\be
F= {\sqrt{ \Big ( 1 - \e^2 \sinh^2(r) \Big ) u^2 
+ \sinh^2(r) v^2 } -\e \sinh^2(r) v \over 1-\e^2 \sinh^2(r) },
\ee
where ${\bf y}= u {\pa \over \pa r}+ v {\pa \over \pa \theta}\in T_{(r, \theta)} ((0, \infty)\times {\rm S}^1)$.
$F$ satisfies ${\bf K}=-1$ and ${\bf S}=0$, but it is not locally projectively flat.

\bigskip
The  Poincare metric on the disk $\Bbb D^2$ is given by 
\be
 \Phi ({\bf y} ) = { 2 \sqrt{u^2+v^2}\over 1-x^2-y^2}, 
\ee
where ${\bf y}= u {\pa \over \pa x}+ v{\pa \over \pa y}\in T_{ (x,y)}\Bbb D^2.$
The Poincare metric has negative constant curvature ${\bf K}=-1$.
Take ${\bf v}=-y {\pa \over \pa x}+ x{\pa \over \pa y}\in T_{(x,y)}\Bbb D^2$ and define $F$ by (\ref{eq}). We obtain
\be
F =
{ \sqrt{\e^2 (xv-yu)^2+(u^2+v^2) \Big ( {1\over 4} (1-x^2-y^2)^2 -\e^2 (x^2+y^2)   \Big )   } -\epsilon ( xv-yu)
\over {1\over 4} (1-x^2-y^2)^2 -\e^2 (x^2+y^2)    }.
\ee 
$F$  satisfies ${\bf K}=-1$ and ${\bf S}=0$, but it is not locally projectively flat.

\bigskip
The Riemannian metric ${\alpha}$ from Theorem \ref{thm1} is given by
\[ \Phi({\bf y}) := \sqrt{ { u^2+v^2 + (xv-yu)^2\over 1+x^2+y^2}  + (xv-yu)^2 },\]
where ${\bf y}= u {\pa \over \pa x} + v {\pa \over \pa y}\in T_{(x, y)} \R^2$. $\Phi$  has constant curvature ${\bf K}=-4$.
Take ${\bf v}_p = -y {\pa \over \pa x}+ x {\pa \over \pa y}$ at $p=(x,y)$
and define $F$ by (\ref{eq}).
We obtain a Randers metric $F=\alpha+\beta$, where 
$\alpha=\alpha({\bf y})$ and $\beta=\beta({\bf y})$ are given by
\begin{eqnarray*}
\alpha : & = & {\sqrt{  u^2+v^2 + (2+x^2+y^2)(xv-yu)^2 -\e^2 (xu+yv)^2 (1+x^2+y^2)} \over \sqrt{1+x^2+y^2} \Big ( 1-\e^2 (x^2+y^2)(1+x^2+y^2) \Big ) }\\
\beta: & = & - {\e (1+x^2+y^2) (xv-yu) \over 1-\e^2 (x^2+y^2)(1+x^2+y^2)}.
\end{eqnarray*}
 $F$ satisfies ${\bf K}=-4$ and ${\bf S}=0$, but it is not locally projectively flat when $\e \not=0$.
}
\end{rem}

\section{Proof of Theorem \ref{thm3}}

Let 
$\Phi$ denote 
the Funk metric on $\Bbb D^2$. It is given by
\[ \Phi ({\bf y}) 
= {\sqrt{ (u^2+v^2) - (xv-yu)^2} + xu+yv \over 1-x^2-y^2},\]
where ${\bf y} = u {\pa \over \pa x} + v {\pa \over \pa y} \in T_{(x, y)} \Bbb D^2$. 
The Finsler metric in Theorem \ref{thm3} is defined by (\ref{FunkFunk}).
Solving the equation (\ref{FunkFunk}), we obtain
\be
F:=  {\sqrt{ u^2+v^2 - \Big ( \e (xu+yv)+ (xv-yu) \Big )^2 }+ (xu+yv) -\e (xv-yu ) \over 1- (1+\e^2)(x^2+y^2) }\\
. \label{Funktypemetric}
\ee
where  ${\bf y}
= u {\pa \over \pa x} + v {\pa \over \pa y}\in T_{(x, y)}\R^2$.
$F =\alpha+\beta $ is a Randers metric on the disk $\Bbb D^2(\rho)$ with 
$\rho =1/\sqrt{1+\e^2}$, where 
 $\alpha$ and $\beta$  are given by
\begin{eqnarray*}
\alpha & = & {\sqrt{ u^2+v^2 - \Big ( \e (xu+yv) + (xv-yu) \Big )^2 }\over 1- (1+\e^2) (x^2+y^2) }\\
\beta & = & { (xu+yv)-\epsilon (xv-yu)\over 1- (1+\e^2)(x^2+y^2) }
\end{eqnarray*}
To compute the curvatures of $F=\alpha+\beta$, we express the Randers metric in a polar coordinate system $x =r \cos\theta, y =r \sin \theta$.
For a vector ${\bf y} = \mu {\pa \over \pa r} +\nu {\pa \over \pa \theta}$, 
$\alpha=\alpha({\bf y})$ and $\beta = \beta({\bf y})$ are given by
\begin{eqnarray*}
\alpha & = & {\sqrt{ \mu^2 + r^2 \nu^2 -r^2 \Big ( r\nu + \e  \mu \Big )^2}\over 1 - (1+\e^2)r^2 }\\
\beta & = & { r \mu -\epsilon r^2 \nu \over 
1 - (1+\e^2) r^2 }.
\end{eqnarray*}
Express $\alpha =\sqrt{ a_{11}\mu^2 + a_{12} \mu \nu + a_{21}\nu \mu + a_{22} \nu^2 }$
and $\beta = b_1 \mu + b_2 \nu$, where
\begin{eqnarray*}
a_{11} & = &  { 1-\e^2 r^2 \over \Big ( 1-(1+\e^2)r^2 \Big )^2},\\
a_{12} & = & -{\e r^3 \over \Big ( 1-(1+\e^2)r^2\Big )^2}\ =\ a_{21},\\
  a_{22} & =&  { r^2 (1-r^2) \over \Big ( 1-(1+\e^2)r^2 \Big )^2} ,
\end{eqnarray*}
\[ b_1 ={r\over 1- (1+\e^2)r^2}, \ \ \ \ \ b_2 = - {\e r^2 \over 1- (1+\e^2) r^2 }.\]
The geodesic coefficients $\bar{G}^1$ and $\bar{G}^2$ of $\alpha$ are given by
\begin{eqnarray*}
\bar{G}^1 & = & { \Big (\e^2 -5 \e^2 r^2 -\e^4 r^2 +2-2r^2\Big ) r \over 
2 \Big (1-(1+\e^2)r^2\Big )^2}\; \mu^2
+ {\e \Big ( 1-r^2+\e^2 r^2 \Big ) r^2\over \Big ( 1-(1+\e^2)r^2 \Big )^2 } \; \mu \nu\\
&& - {\Big (1-r^2+\e^2 r^2\Big )(1-r^2) r \over 2\Big (1-(1+\e^2)r^2\Big )^2 }\; \nu^2\\
\bar{G}^2 & = & {  \e \Big ( -3 +r^2+3\e^2 r^2 \Big )   \over 2\Big (1-(1+\e^2)r^2\Big )^2 }\; \mu^2
+ { (1-\e^2r^2)\Big (1-r^2+\e^2r^2 \Big ) \over \Big ( 1- (1+\e^2)r^2 \Big )^2 r } \; \mu \nu \\
&&
 - {\e \Big ( 1-r^2 +\e^2 r^2 \Big ) r^2 \over 
2 \Big ( 1- (1+\e^2)r^2 \Big )^2 }\; \nu^2
\end{eqnarray*}
The Gauss curvature $\bar{\bf K}$ of $\alpha$ is given by
\be
 \bar{\bf K} = { -1 -5 \e^2 + (1-\e^4) r^2 \over 1 - (1+\e^2) r^2 }.\label{K2}
\ee
We see that $\bar{\bf K}$ is not a constant unless $\e =0$.

Now we are going to find the geodesic coefficients $G^1$ and $G^2$ of $F=\alpha+\beta$. Let $r_{ij}, s_{ij}, s^i_{\ j}, s_j $ and $e_{ij}$ as above. 
A direct computation yields that 
\begin{eqnarray*}
r_{11} & = & { 1-r^2 - 3\e^2 r^2 \over 
\Big ( 1-(1+\e^2)r^2 \Big )^2 }\\
r_{12} & = & - {\e (1-\e^2) r^3\over \Big (1-(1+\e^2)r^2\Big )^2} \ \ = \ \ r_{21}\\
r_{22} & = & { \Big ( 1-r^2 +\e^2 r^2 \Big ) r^2 
\over \Big ( 1-(1+\e^2)r^2 \Big )^2}\\
s_{11} & = & 0 \ \ =\ \  s_{22}\\
s_{12} & = & {\e r \over \Big ( 1 - (1+\e^2)r^2 \Big )^2 } \ \ =  \ \ -s_{21}\\
s^1_{\ 1} & = &  -{\e^2r^2 \over 1- (1+\e^2) r^2 }\ \ = \ \ - s^2_{ \ 2} \\
s^1_{\ 2} & = & { \e r (1-r^2) \over 1- (1+\e^2)r^2 }\\
s^2_{\ 1} & = & - {\e ( 1-\e^2 r^2 ) \over (1- (1+\e^2) r^2) r }\\
s_1 & = & {\e^2 r \over 1- (1+\e^2) r^2 }\\
s_2 & = & {\e r^2 \over 1- (1+\e^2) r^2 }.
\end{eqnarray*}
We immediately see that 
\be
 e_{ij} :=r_{ij}+b_i s_j +b_j s_i= a_{ij}- b_i b_j.\label{eij}
\ee
By Lemma 3.1 in \cite{ChSh}, 
(\ref{eij}) is equivalent to that
\[ {\bf S}={3\over 2} F.\]
 By (\ref{Gi2}) and  the above identities, we obtain
\begin{eqnarray*}
G^1 & = & \bar{G}^1+{1\over 2} (\alpha -\beta)\mu - {\e r \Big ( \e \mu + r \nu \Big ) 
\over 1 - (1+\e^2) r^2 } \; \mu -  {\e r \Big ( \e r \mu -\nu + r^2 \nu \Big ) 
\over 1 - (1+\e^2) r^2 } \; \alpha \\
G^2 & = & \bar{G}^2 + {1\over 2} (\alpha -\beta)\nu -{\e r \Big ( \e \mu + r \nu \Big ) 
\over 1 - (1+\e^2) r^2 }\; \nu  - { \e \Big ( \mu -\e^2 r^2 \mu - \e r^3 \nu \Big ) \over r \Big ( 1- (1+\e^2) r^2 \Big ) }\; \alpha  
\end{eqnarray*}
Plugging them into (\ref{Rik}), we immediately obtain
\[ R^i_{\ k} = -{1\over 4} \Big \{ F^2 \delta^i_k - F F_{y^k} y^i\Big \}.\]
Thus the Gauss curvature ${\bf K}=-1/4$.

We can also use (\ref{RikRik}) to verify that ${\bf K}= -1/4$. 
To do so, it suffices to compute $s_{0|0}$ and $s^k_{\ 0|k}$. They are given by
\begin{eqnarray*}
s_{0|0} & = & { \e^2 ( 1 + r^2 -\e^2 r^2 )\over \Big ( 1- (1+\e^2) r^2\Big  )^3    }\; \mu^2 - { 4 \e^3 r^3 \over \Big ( 1- (1+\e^2) r^2 \Big )^3 }\; \mu\nu 
+ { \e^2 r^2 ( 1- r^2 + \e^2 r^2)   \over \Big ( 1- (1+\e^2) r^2 \Big )^3 }\; \nu^2  \\
s^k_{\ 0|k} & = & - { \e^2 (1+\e^2)  r^3 \over\Big  ( 1- (1+\e^2) r^2\Big  )^2 }\; \mu 
+ {\e (1+\e^2)  r^2 ( 1 - r^2  ) \over\Big  ( 1- (1+\e^2) r^2\Big  )^2 }\; \nu.
\end{eqnarray*} 
Plugging $c= 1/2$ and the above identities into (\ref{RikRik}) gives 
\[ {\bf Ric} = - {1\over 4} F^2.\]
We conclude that ${\bf K} = {\bf Ric}/F^2 = -1/4$.

\bigskip

\begin{rem}
{\rm Below is a byproduct.
Let
\begin{eqnarray*}
\alpha : & = & {\sqrt{ u^2+v^2 - \Big ( \e (xu+yv) + (xv-yu) \Big )^2 }\over 1- (1+\e^2) (x^2+y^2) }\\
\tilde{\alpha} : & = & {\sqrt{ \Big(1-(1+\e^2)(x^2+y^2)\Big )(u^2+v^2)
-\Big ( 1-\e^2 -(x^2+y^2) \Big ) (xv-yu)^2 }
\over \Big ( 1- (1+\e^2)(x^2+y^2) \Big ) \sqrt{ 1-x^2-y^2} }.\end{eqnarray*}
$\alpha$ and $\tilde{\alpha}$ are two Riemannian metrics
on $\Bbb D^2(\rho)$ with radius $\rho = 1/\sqrt{1+\e^2}$.
According to (\ref{K1}) and (\ref{K2}), 
The Gauss curvatures  of $\alpha$ and $\tilde{\alpha}$
are equal and given by
\[ \bar{\bf K}=
{ -1 -5 \e^2 + (1-\e^4) (x^2+y^2)
\over 1- (1+\e^2) (x^2+y^2) }.
\]
}
\end{rem}

\noindent
Math Dept, IUPUI, 402 N. Blackford Street, Indianapolis, IN 46202-3216, USA.  \\
zshen@math.iupui.edu

\end{document}